\documentstyle{article}
\newtheorem{df}{ \sc Definition}[section]

\newtheorem{pr}[df]{ \sc Proposition}

\newtheorem{th}[df]{ \sc Theorem}

\newtheorem{re}[df]{ \it Remark}

\newtheorem{eq}[df]{\rm (}
\def\codim{{\rm codim}\,}

\def\mpr#1{\;\smash{\mathop{\hbox to 20pt{\rightarrowfill}}\limits^{#1}}\;}
\def\epi#1{\;\smash{\mathop{\hbox to 20pt{\rightarrowfill}\hskip
-13pt\rightarrow}\limits^{#1\,}}\;}
\def\epii{\smash{\mathop{\hbox to 14pt{\rightarrowfill}\hskip
-10pt\rightarrow}}}
\def\mono{\lhook\joinrel\relbar\joinrel\rightarrow}
\def\mpL#1{\;\smash{\mathop{\hbox to 20pt{\leftarrowfill}}\limits_{#1}}\;}
\def\mpl#1{\;\smash{\mathop{\hbox to 20pt{\leftarrowfill}}\limits^{#1}}\;}

\def\Z{{\bf Z}}

\def\Q{{\bf Q}}
\def\C{{\bf C}}

\def\F{{\bf F}}
\def\C{{\bf C}}
\def\Tor{{\rm Tor}}

\def\s{\vskip6pt}
\def\x{{*}}
\def\Proof{\noindent{\it Proof. }}
\def\Proo{\s\noindent{\it Proof~}}
\def\qed{\hfill$\Box$\vskip10pt}

\begin{document}

\author{Andrzej Weber\footnote{Supported by KBN grant  1 P03A 005 26 }
\\
\it Instytut Matematyki, Uniwersytet Warszawski}

\title{ON TORSION IN HOMOLOGY OF SINGULAR TORIC VARIETIES}

\date{August 2005}
\maketitle

\begin{abstract} Let $X$ be a toric variety. Rationally
Borel-Moore homology of $X$ is isomorphic to the homology of the
Koszul complex $A^T_*(X)\otimes \Lambda^\x M$, where $A^T_*(X)$ is
the equivariant Chow group and $M$ is the character group of $T$.
Moreover, the same holds for coefficients which are
the integers with certain primes inverted.
\end{abstract}

\section{Introduction}
Let $G$ be a connected Lie group which acts on a topological space
$X$. The equivariant cohomology of $X$  is defined to be the
cohomology of the space obtained from Borel construction
 $$H^*_G(X)=H^*(EG\times_GX)\,.$$
It often happens that this group is easily computable. This is so
for example when $X$ is a smooth complex algebraic manifold acted
by a complex algebraic group with only finitely many orbits. Then
the equivariant cohomology with rational coefficients is the
direct sum

\begin{eq}\rm ) $$H^*_G(X;\Q)=\bigoplus_{\rm orbit}
H^{*-2c}(BG_x;\Q)\,.$$\end{eq}
 Here $G_x$ is the stabilizer of a
point from the orbit and $c$ is the complex codimension of the
orbit.

If $X$ is singular then there is another invariant which is easy
to compute. That is the equivariant Borel-Moore homology,
\cite[\S2.8]{EG}. It can be interpreted as the equivariant
cohomology with coefficients in the dualizing sheaf. The
equivariant Borel-Moore homology usually is nontrivial in the
negative degrees.
 Again we have
 \begin{eq}\rm ) \label{row} $$H^{BM,G}_*(X;\Q)=\bigoplus_{\rm orbit}
H^{2d-*}(BG_x;\Q)\,.$$ \end{eq}
 Here $d$ is the complex dimension
of the orbit. This formula, as well as the previous one, follows
from the fact that the rational cohomology of the classifying
space is concentrated in even degrees.

A passage from the equivariant cohomology to the usual one is
possible due to the Eilenberg-Moore spectral sequence. The second
table is of the form
 $$E^{-p,q}_2=
\Tor^{q,H^*(BG;\Q)}_p(H^*_G(X;\Q),\Q)\Rightarrow
H^{q-p}(X;\Q)\,.$$
 (Torsion functor has two gradings: $p$ is the usual grading of the left derived
functor and $q$ is the internal grading.) A generalization of this
spectral sequence for cohomology with sheaf coefficients was
described in \cite{FW}. If $G$ is a complex algebraic group, $X$
is an algebraic variety and the action is algebraic then all the
cohomological invariants are equipped with the weight filtration.
This often forces the spectral sequence to degenerate. In
particular \cite[Th.~1.6]{FW} we have:
\begin{th} If $X$ is smooth and the action has only finitely many
orbits then the rational cohomology of $X$ is given additively by:
$$H^i (X;\Q)=\bigoplus_{q-p=i} \Tor^{q,H^*(BG;\Q)}_p(H^*_G(X;\Q),\Q)\,.$$
\end{th}
Having the decomposition \ref{row}  and using the fact that the
equivariant cohomology of an orbit is pure we apply
\cite[Th.~1.3]{FW}. We obtain
\begin{th} If the action has only finitely many
orbits then the rational Borel-Moore homology of $X$ is given by:
$$H_i^{BM} (X;\Q)=\bigoplus_{p-q=i} \Tor^{q,H^*(BG;\Q)}_p(H^{BM,G}_*(X;\Q),\Q)\,.$$
\end{th}
In both theorems above we assume that spaces and actions are
algebraic.

In this note we want to specialize our results to toric varieties.
The intersection cohomology was already described in \cite{W}.
Here we study ordinary homology but, in addition, we care about
the integral coefficients. First we note that for a torus $T$ the
Eilenberg-Moore spectral sequence can be replaced by an easier one
(which is in fact isomorphic after a renumbering of entries). This
is just the spectral sequence of the fibration
$$T\subset ET\times X \epi{} ET\times_T X\,.$$
The second table
 $$E^{p,q}_2=H^p_T(X)\otimes H^q(T)$$
 with its differential is exactly the Koszul complex. Therefore
 $$E_3^{p,q}=\Tor^{p+2q,H^*(BG)}_q(H^*_G(X),\Z).$$
The exact degrees are slightly surprising, but they agree with the
weight filtration when $X$ is smooth. For singular varieties we
will describe the homological variant of the spectral sequence in
an elementary way.

Now we apply the Frobenius  endomorphism of the toric variety.
This allows to show that the spectral sequence degenerates not
only over $\Q$ but also with small primes inverted. We prove two
theorems.

\begin{th} \label{teo1} The
above spectral sequence degenerates on $E^3$ for rational
coefficients and for the coefficients in $\F_q$ if $q>\lceil{\dim
X \over 2}\rceil$.\end{th}

\begin{th}\label{teo2} If $q>\lceil{\dim X+1\over 2}\rceil$ then the
$q$-torsion of the integral homology is a direct sum of the
$q$-torsions in $E^3$.\end{th}

As a consequence we obtain:

\begin{th} Let $X$ be a toric variety and let $R$ be the ring of integers
with inverted
primes which are smaller or equal to  $\lceil{\dim X+1\over
2}\rceil$.
Then
 $$H_i^{BM} (X;R)=\bigoplus_{p-q=i} \Tor^{q,H^*(BG;R)}_p(H^{BM,T}_*(X;R),R)\,.$$
\end{th}

It remains to remark that $H_{2i}^{BM} (X;R)\simeq A_i^T(X)\otimes
R$ is the equivariant Chow group.

We suspect that the assumption about $q$ is redundant.

\section{Equivariant Borel-Moore homology}

From now on we omit coefficients in the notation. Let $X$ be an
algebraic variety acted by the torus $T=(\C^*)^n$. Let
 $$ET_d=(\C^{d+1}\setminus \{0\})^n$$
 be an approximation of $ET$. The equivariant Borel-Moore homology
is defined by the formula:
$$H_i^{BM,T}(X)=\lim_{\stackrel{\longleftarrow}{d}} H_{i+2nd}^{BM}(ET_d\times_T X)\,.$$
The limit is taken with respect to the inverse system
$$\iota_d^!:H_{i+2n(d+1)}^{BM}(ET_{d+1}\times_T X)\to
H_{i+2nd}^{BM}(ET_d\times_T X)\,.$$ The Gysin map $\iota_d^!$ is
defined since the inclusion (of the real codimension $2n$)
$$\iota_d:ET_{d}\times_T X\mono ET_{d+1}\times_T X$$ is normally nonsingular.
In fact the limit stabilizes:
 $H_i^{BM,T}(X)=H_{i+2d}^{BM}(E_{d+1}\times_T X)$ for
$i>(1-2d)+2\dim X$.

Let $X$ be a toric variety. By \ref{row} or \cite{BZ} the equivariant
homology
is the sum of
homologies of orbits:
\begin{eq}\rm )\label{row2}\hfil
$$H_i^{BM,T}(X)=\bigoplus_{\sigma\in
\Sigma}H_i^{BM,T}(O_\sigma)\,.$$
\end{eq}

 \noindent The orbits $\sigma$ are labelled by
a fan $\Sigma$. The equivariant Borel-Moore homology is isomorphic
to the equivariant Chow group $A^T_*(X)$ considered e.g.~by Brion
\cite{Br}. Each $H_*^{BM,T}(O_\sigma)$ is isomorphic to
$Sym(\langle\sigma\rangle)[-2\codim \sigma]$ (the symmetric power
$Sym^i(\langle\sigma\rangle)$ is placed in the degree
$2(\codim\sigma-i)$). In particular it is a free abelian group.
The odd part of the equivariant Borel-Moore homology vanishes. The
modul structure over $H^*(BT)$ is described in \cite{Br}.

\section{Frobenius endomorphism}
Let $p>1$ be a natural number. Toric varieties are equipped with
Frobenius endomorphism (power map) $\phi_p:X\to X$. The power map
of $T$ is denoted by $\psi_p$. Both maps induce a map at each step
of the approximation of the Borel construction $ET_d\times_TX$. We
denote this map by $\phi_p^{T}$. (We note that $\phi^T_p$ is not
the same as $1\otimes \phi_p$ considered in \cite{BZ}.)

We would like to encode the action of $\phi_p^T$ in equivariant
homology. The map $\iota_d^!$ does not commute with
$\phi_{p,*}^{T}$, but
$\phi_{p,*}^{T}\iota_d^!=p\,\iota_d^!\phi_{p,*}^{T}$. Then we set
$H_i^{BM,T}(X)= H_{i+2d}^{BM}(ET_d\times_T X)(d)$ for $d$
sufficiently large. If we study homology with coefficients in the
field $\F=\Q$ or $\F_q$, provided that $(p,q)=1$, the symbol $(d)$
denotes tensoring with $\F$ acted by $\phi_p$ via the
multiplication by $p^{-d}$. If we want to study integral homology
we just analyze homology of a sufficiently large approximation of
the Borel construction. Now we can state
\begin{pr} $\phi_{p,*}^{T}$ induces the multiplication by $p^i$ on
$H_{2i}^{BM,T}(X)$. \end{pr}

\Proof Homology of an orbit is generated over $H^*(BT)$ by its
fundamental class. The map $\phi_{p,*}^T$ restricted to the orbit
$\cal O_\sigma$ is a covering of the degree $p^i$, where
$i=\dim{\cal O_\sigma}=\codim\,\sigma$. Thus
$\phi_{p,*}^T[\sigma]=p^i[\sigma]\in H_{2i}^{BM,T}({\cal
O}_\sigma)$. Now we apply the additivity $\ref{row2}$. \qed

We consider a system of $T$-\-fibrations
$$ET_d\times X\to ET_d\times_T X\,.$$
We obtain a system of spectral sequences with
$${_dE}^2_{k,l}=H^{BM}_k(ET_d\times_T X)\otimes H^{BM}_l(T)\,.$$ The
map $\iota_d^!$ passes to a map of spectral sequences
 $$_{d+1}E_{k+2,l}^j(1)\to{_dE}_{k,l}^j\,.$$
We set ${_\infty E}_{k,l}^j={_dE}_{k+2nd,l}^j(d)$ for $d$ large
enough. The map $\phi_p^*$ of $H^{BM}_q(T)$ is the multiplication
by $p^{q-n}$ (for $q<n$ this group is trivial). Therefore the
resulting map of the spectral sequences is the multiplication by
$p^{{k\over 2}+l-n}$ on ${_\infty E}_{k,l}^2$.

\begin{th} \label{teo} The spectral sequence ${_\infty E}_{k,l}^j$ converges
to
$H^{BM}_{k+l-2n}(X)$. For rational coefficients it degenerates on
$E^3$. The resulting filtration
coincides after renumbering with the weight filtration of
homology.\end{th}

\begin{re} \rm There is a shift $-2n$ in the degree which is
repaired when we move the generators of $H^{BM}_*(T)$ to the
negative degrees. They should be placed there since we compare
$E_2$ with the Koszul complex, see below.\end{re}

\begin{re} \rm The homology $H^{BM,T}_*(X)$ is a module over
$H^*(BT)=Sym\, M$, where $M=Hom(T,\C^*)$. The differential
$d_2:{_\infty E}_{k,l}^2\to {_\infty E}_{k-2,l+1}^2$ after the
identification $H_*^{BM}(T)=H^{2n-*}(T)$ with $\Lambda^\x M$
becomes the Koszul differential
$$H^{BM,T}_k(X)\otimes \Lambda^{2n-l}M\to
H^{BM,T}_{k-2}(X)\otimes \Lambda^{2n-l-1}M
\,.$$
\end{re}

\begin{re}\rm Koszul complex contains a complex constructed
in \cite{T}. Totaro considered rational coefficients, but he
remarked that some information about the torsion can be obtained.
\end{re}

\begin{re}\rm From another point of view the Koszul duality and
Frobenius endomorphism appears in \cite{B} and \cite{BL} for dual
pairs of affine toric varieties.\end{re}

\Proo {\it of Theorem \ref{teo}}. At each step ${_d}E^*_{k,l}$
converges to $H^{BM}_{k+l}(X\times(\C^{d+1}\setminus\{0\})^n)$,
which is equal to
$H^{BM}_{k+l}(X\times(\C^{(d+1)n})=H^{BM}_{k+l-2n(d+1)}(X)$ for
$d$ sufficiently large. Therefore ${_\infty E}^*_{k,l}$, which is
equal to ${_d E}^*_{k+2nd,l}$ (for $d$ sufficiently large)
converges to $H^{BM}_{k+l-2n}(X)$. The eigenvalue of $\phi_p$
acting on ${_\infty E}^*_{k-r,l+r-1}$ is equal to $p^{{k\over 2} +
l -n + {r\over 2}-1}$. There is no obstruction for the
differential $d_2$, but the higher differentials have to
vanish.\qed

\section{Torsion.}

It is not possible to detect the $q$-torsion of $X$ for small
prime $q$, but we prove Theorems \ref{teo1} and \ref{teo2}
announced in the introduction:

\Proo {\it of Theorem \ref{teo1}}.  To show that the higher
differentials vanish we consider the eigenvalues of $\phi_p$
acting on ${_dE}^2_{k-r,l+r-1}$ for $r=2,\dots n+1$ and $k-r$
even. These are at most $\lceil {n\over 2}\rceil$ subsequent
powers of $p$, as the reader may easily check (see the picture
below). We will chose $p$ such that all these powers are different
modulo $q$. It is enough to take $p$ which generate the group
$\F^*_q\simeq \Z/(q-1)$. \qed

\vbox{
$$\matrix{2n&\dots&\circ&\cdot&\bullet&\cdot&\circ&\cdot&\circ&\cdot&\circ&\cdot&\circ&\dots\cr
      &\dots&\circ&\cdot&\circ &\cdot&\circ&\cdot&\circ&\cdot&\circ&\cdot&\circ&\dots\cr
      &\dots&\circ&\cdot&\circ &\cdot&\bullet&\cdot&\circ&\cdot&\circ&\cdot&\circ&\dots\cr
      &\dots&\circ&\cdot&\circ &\cdot&\circ&\cdot&\circ&\cdot&\circ&\cdot&\circ&\dots\cr
      &\dots&\circ&\cdot&\circ &\cdot&\circ&\cdot&\star&\cdot&\circ&\cdot&\circ&\dots\cr
      &\dots&\circ&\cdot&\circ &\cdot&\circ&\cdot&\circ&\cdot&\star&\cdot&\circ&\dots\cr
     n &\dots&\circ&\cdot&\circ &\cdot&\circ&\cdot&\circ&\cdot&\circ&\cdot&\circ&\dots\cr}$$
\s
 \hfil
 \hbox{\footnotesize\vbox{
 \centerline{ A picture of the spectral sequence for $n=6$.}
 \centerline{$\star$ denotes the source and the target of $d_2$ (which always
 have the same weight),}
 \centerline{$\bullet$ denotes the remaining possibly
 nonzero entries of the spectral sequence which}
 \centerline{are hit by the
higher differentials.}}} }\s

If we want to determine the $q$-torsion for integral homology we
may meet problems with extensions.

\Proo {\it of Theorem \ref{teo2}}. The part of the integral
spectral sequence computing the $q$-torsion degenerates as in the
previous proof. To avoid problems with extensions we have to know
that the eigenvalues of $\phi_p$ on $E^3_{k,l}$ are different
along the lines $k+l=const$. There are exactly $\lceil{n+1\over
2}\rceil$ subsequent powers of $p$. We proceed as before, that is
we find $p$ with different powers modulo $q$. \qed

\vbox{
$$\matrix{2n&\dots&\circ &\cdot&\bullet&\cdot&\circ&\cdot&\circ&\cdot&\circ&\cdot&\circ&\dots\cr
      &\dots&\circ &\cdot&\circ &\cdot&\circ&\cdot&\circ&\cdot&\circ&\cdot&\circ&\dots\cr
      &\dots&\circ &\cdot&\circ &\cdot&\bullet&\cdot&\circ&\cdot&\circ&\cdot&\circ&\dots\cr
      &\dots&\circ &\cdot&\circ &\cdot&\circ&\cdot&\circ&\cdot&\circ&\cdot&\circ&\dots\cr
      &\dots&\circ &\cdot&\circ &\cdot&\circ&\cdot&\bullet&\cdot&\circ&\cdot&\circ&\dots\cr
      &\dots&\circ &\cdot&\circ &\cdot&\circ&\cdot&\circ&\cdot&\circ&\cdot&\circ&\dots\cr
     n &\dots&\circ &\cdot&\circ &\cdot&\circ&\cdot&\circ&\cdot&\bullet&\cdot&\circ&\dots\cr}$$
\s
 \hfil
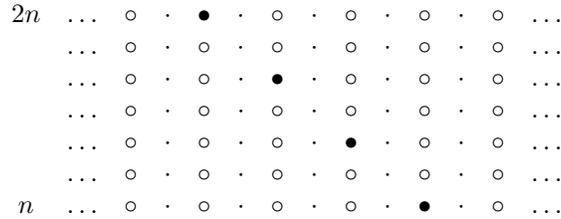
 \hbox{\footnotesize\vbox{\centerline{A picture of the spectral sequence
for $n=6$.}
               \centerline{ The entries $\bullet$ should have different weights
                          in $\F_q$.}}} }
 \s

We conjecture that the theorems above are true without assumptions
on $q$. The conjecture holds if $X$ is smooth by the work of
M.~Franz \cite{F}.

\section{Weight filtration and gradation.} Let $V_j$ be the
eigensubspace of $\phi_p$ acting on $H^{BM}_*(X;\Q)$ for the
eigenvalue $p^j$. It does not depend on $p$. The weight filtration
in homology usually is denoted by $W^i=W_{-i}$ and it is
decreasing. Our gradation is related to the weight filtration:
$$W^i=\bigoplus_{j\geq {i\over 2}} V_j\,.$$
In particular $W^{2j}=W^{2j-1}$.

\s\s{\noindent\sc\hskip4pt Conjugation.} Toric varieties are
defined over real numbers. The complex conjugation acts on the
complex  points  of  $X$.  We  can  also determine the action on
the homology: it acts by $(-1)^{i\over 2}$ on the $i$-th term of
the weight gradation.

\s
{\parindent 140pt

 Andrzej Weber

Instytut Matematyki, Uniwersytet Warszawski,

ul. Banacha 2, 02-097 Warszawa, POLAND

e-mail: \tt aweber@mimuw.edu.pl}

\begin{thebibliography}{BBB}

\bibitem [B]{B} T.~Braden {\em Koszul duality for toric varieties },
 arXiv.org:math/0308216


\bibitem [BL]{BL} T.~Braden, V.~A.~Lunts {\em
 Equivariant-constructible Koszul duality for dual toric
 varieties},
arXiv.org:math.AG/0409495

\bibitem [Br]{Br} M.~Brion, {\em Equivariant Chow groups for torus
actions} Transformation groups, Vol 2, No 3 (1997), pp. 1-43


\bibitem [BZ]{BZ} J-L.~Brylinski, B.~Zhang: {\em Equivariant Todd
classes for toric varieties}, arXiv.org:math.AG/0311318


\bibitem[EG]{EG} D.~Edidin, W.~Graham {\em Equivariant
intersection theory} Invent. Math. vol 131 (1998) no. 3, pp.
595--634

\bibitem [F]{F} M.~Franz {\em
On the integral cohomology of smooth toric varieties},
arXiv.org:math.AT/0308253

\bibitem [FW]{FW}
M.~Franz, A.~Weber {\em Weights in cohomology and the
Eilenberg-Moore spectral sequence},   Ann. Inst Fourier (Grenoble)
55 (2005), no. 2, 673--691


\bibitem [T]{T} B.~Totaro {\em Chow groups, Chow cohomology, and
linear varieties}. Journal of Algebraic Geometry, to appear.\\
http://www.dpmms.cam.ac.uk/$\sim$bt219/papers.html


\bibitem [W]{W} A.~Weber {\em
Weights in the cohomology of toric varieties}.  Cent. Eur. J.
Math. 2 (2004), no. 3, 478--492, arXiv.org:math.AG/0301314

\end{thebibliography}
\end{document}